\documentclass[11pt]{amsart} 

\usepackage[lmargin=1in,rmargin=1in,tmargin=1in,bmargin=1in]{geometry}
\usepackage[ps,all,arc,rotate]{xy}
\usepackage{graphicx, float, epstopdf}
\usepackage{hyperref, color}
\usepackage{centernot}
\usepackage{fancyhdr}
\usepackage[utf8]{inputenc}
\usepackage{amsfonts,amssymb,amsmath,amsthm,mathrsfs}
\usepackage{graphics, setspace}
\usepackage{braket}
\usepackage{mathtools}
\usepackage{tikz}

\numberwithin{equation}{section}
\numberwithin{figure}{section}      

\newtheorem{lemma}{Lemma}[section]
\newtheorem{theorem}{Theorem}[section]

\newtheorem{conjecture}[lemma]{Conjecture}
\newtheorem{corollary}{Corollary}[section]
\theoremstyle{definition}

\newtheorem{remark}{Remark}[section]

\newcommand{\bZ}{\mathbb{Z}}

\newcommand{\bQ}{\mathbb{Q}}
\newcommand{\bA}{\mathbb{A}}

\newcommand{\cF}{\mathcal{F}}
\newcommand{\ve}{\varepsilon}

\newcommand{\bR}{\mathbb{R}}
\newcommand{\bC}{\mathbb{C}}

\newcommand{\cC}{\mathcal{C}}

\newcommand{\cS}{\mathcal{S}}

\newcommand{\cT}{\mathcal{T}}
\newcommand{\Tr}{\text{Tr}}
\newcommand{\cJ}{\mathcal{J}}

\newcommand{\bH}{\mathbb{H}}

\usepackage{hyperref}
\let\svthefootnote\thefootnote
\newcommand\freefootnote[1]{%
  \let\thefootnote\relax%
  \footnotetext{#1}%
  \let\thefootnote\svthefootnote%
}

\usepackage[
backend=biber,
style=alphabetic,
]{biblatex}
\title{Weighted Low-lying Zeros of L-functions Attached to Siegel Modular Forms}
\author{Shifan Zhao}

\begin{document}

\freefootnote{2020 Mathematics Subject Classification: Primary 11F46, 11F66, 11F72  }

\maketitle

\begin{abstract}
In this paper, we study weighted low-lying zeros of spinor and standard $L$-functions attached to degree 2 Siegel modular forms. We show the symmetry type of weighted low-lying zeros of spinor $L$-functions is symplectic, for test functions whose Fourier transform have support in $(-1,1)$, extending the previous range $(-\frac{4}{15},\frac{4}{15})$ in \cite{KowalksiSahaTsimerman2012}. We then show the symmetry type of weighted low-lying zeros of standard $L$-functions is also symplectic. We further extend the range of support by performing an average over weight. As an application, we discuss non-vanishing of central values of those $L$-functions.
\end{abstract}

\tableofcontents

\section{Introduction}

D. Hilbert and G. P\'olya suggested that non-trivial zeros of the Riemann zeta function $\zeta(s)$ correspond to eigenvalues of a self-adjoint operator on some Hilbert space. The first evidence of such a connection was found by H. L. Montgomery \cite{Montgomery1973}, who investigated the pair correlation of non-trivial zeros of $\zeta(s)$ and conjectured that it is, as pointed out by F. J. Dyson, the same as the pair correlation of eigenvalues of random Hermitian or unitary matrices of large order, also known as the Gaussian Unitary Ensemble (GUE) model. This conjecture of Montgomery was later supported by numerical results by A. M. Odlyzko \cite{Odlyzko1987}, based on values for the first $10^5$ zeros and for zeros number $10^{12}+1$ to $10^{12}+10^5$. The local spacing between these sample zeros matches prediction by the GUE model quite well.

Z. Rudnick and P. Sarnak \cite{RudnickSarnak1996} extended Montgomery's work by computing the general $n$-level correlation function of zeros of any principal $L$-function $L(s,\pi)$ attached to a cuspidal automorphic representation $\pi$ of GL$_m(\bA_\bQ)$ (in a restricted range). Their answer is universal and is precisely the one predicted by the GUE model. Numerical evidence were found by R. Rumely \cite{Rumely1993} for primitive Dirichlet $L$-functions, and by M. O. Rubinstein \cite{Rubinstein1998} for Hasse-Weil $L$-functions of 3 distinct elliptic curves and for the Hecke $L$-function associated to Ramanujan's $\tau$-function.

Although the $n$-level correlation statistic of zeros of any fixed principal automorphic $L$-function obeys the universal GUE law, there is another statistic, called the $n$-level density of low-lying zeros, that is sensitive to families. N. Katz and Sarnak \cite{KatzSarnak1999} studied low-lying zeros of zeta functions of varieties over finite fields (the "function field" analogue). For these they indicated that a spectral interpretation exists in terms of eigenvalues of Frobenius on cohomology groups. On the number field side, although many results concerning low-lying zeros have been proved, it is still not clear where their spectral nature comes from. See also \cite{KatzSarnak19992} for a nice survey on these topics. 

Before stating our results, we first describe the problem in general terms. Let $\cF_Q$ be a family of automorphic forms, ordered by certain conductor $Q \geq 1$. To each $f \in \cF_Q$ one associates an $L$-function
\begin{equation}
    L(s,f) = \sum_{n=1}^\infty \frac{\lambda_f(n)}{n^s},
\end{equation}
which converges absolutely for $s \in \bC$ in certain right half-plane. We assume that $L(s,f)$ continues analytically to the whole complex plane $\bC$ with possibly finitely many poles and satisfies a functional equation
\begin{equation}
    \Lambda(s,f) = L_\infty(s,f)L(s,f) = \ve_f\Lambda(1-s,f),
\end{equation}
where $\ve_f = \pm 1$ is the root number.

We also assume Generalized Riemann Hypothesis (GRH) for $L(s,f)$. That is, non-trivial zeros of $L(s,f)$ all lie on the critical line. Thus we may denote those zeros by
\begin{equation}
    \rho_f = \frac{1}{2} + i\gamma_f, \hspace{3mm} \gamma_f \in \bR.
\end{equation}

Let $\Phi \in \cS(\bR)$ be an even Schwartz function whose Fourier transform $\hat{\Phi}$ has compact support. We call such a function a "test function" throughout this paper. To this end we define the 1-level density of low-lying zeros of $L(s,f)$, with respect to the test function $\Phi$, to be
\begin{equation}
    D(f;\Phi) = \sum_{\rho_f} \Phi\left(\frac{\gamma_f}{2\pi}\log c_f\right),
\end{equation}
where $\rho_f$ runs through non-trivial zeros of $L(s,f)$, counted with multiplicity, and $c_f$ is a parameter associated with $f \in \cF_Q$, comparable to the analytic conductor of $f$ (to be specified later). The Density Conjecture for low-lying zeros of $L(s,f)$ asserts the following:
\begin{conjecture}[Density Conjecture]\label{density conjecture}
    For any even Schwartz function $\Phi$ whose Fourier transform $\hat{\Phi}$ has compact support, we have 
    \begin{equation}
    \lim_{Q \to \infty} \frac{1}{|\cF_Q|}\sum_{f \in \cF_Q} D(f;\Phi) = \int_{-\infty}^\infty \Phi(x)W(\cF)(x)dx  
    \end{equation}
    for some distribution $W(\cF)$ depending only on $\cF$.
\end{conjecture}

Many observations and results in \cite{KatzSarnak1999} suggest that the distribution $W(\cF)$ depends on the family $\cF$ through a symmetry group $G(\cF)$. Those symmetry types are orthogonal O, special orthogonal even SO(even), special orthogonal odd SO(odd), symplectic Sp and unitary U, with corresponding distributions and Fourier transforms
\begin{align}
    W(O)(x) = 1+\frac{1}{2}\delta_0(x)&, \hspace{3mm} \hat{W}(O)(y) = \delta_0(y) + \frac{1}{2}, \\
    W(SO(\text{even}))(x) = 1+ \frac{\sin 2\pi x}{2\pi x}&, \hspace{3mm} \hat{W}(SO(\text{even}))(y) = \delta_0(y) + \frac{1}{2}\eta(y),  \\
    W(SO(\text{odd}))(x) = 1-\frac{\sin 2\pi x}{2\pi x} + \delta_0(x)&, \hspace{3mm} \hat{W}(SO(\text{odd}))(y) = \delta_0(y) - \frac{1}{2}\eta(y) + 1, \\
    W(Sp)(x) = 1 - \frac{\sin 2\pi x}{2\pi x}&, \hspace{3mm} \hat{W}(Sp)(y) = \delta_0(y) - \frac{1}{2}\eta(y), \label{Fourier pair} \\
    W(U)(x) = 1&, \hspace{3mm} \hat{W}(U)(y) = \delta_0(y),
\end{align}
where $\delta_0$ is the Dirac distribution at $0$, and $\eta(y) = 1,\frac{1}{2},0$ for $|y|<1, |y|=1$ and $|y|>1$ respectively. The first three distributions of different orthogonal symmetry type have indistinguishable Fourier transforms within $(-1,1)$, while the symplectic and unitary symmetry types are distinguishable from the orthogonal ones.  

The Density Conjecture \ref{density conjecture} has been verified for many families (in restrict ranges). See \cite{IwaniecLuoSarnak2000,Rubinstein2001,FouvryIwaniec2003,Guloglu2005,Young2006,DuenezMiller2006,DuenezMiller2006,GaoZhao2011,ChoKim2015,ShinTemplier2016,LiuMiller2017,KimWakatsukiYamauchi2020}, to name a few. In almost all results towards this direction, the support of Fourier transform of the test function $\Phi$ is restricted within certain range. One important question in this topic is how to extend the range as large as possible, for the full Density Conjecture \ref{density conjecture} does not require any condition on the compact support of $\hat{\Phi}$.

One can also consider "weighted" distribution of low-lying zeros by allowing certain weight $\omega_f$. The weighted average density under consideration is
\begin{equation}
    \left(\sum_{f \in \cF_Q}\omega_f\right)^{-1} \sum_{f \in \cF_Q} \omega_f D(f;\Phi).
\end{equation}
Often these weights $\omega_f$ contain important arithmetic information such as central values of $L$-functions, and including them may possibly change the symmetry type. Recent results in this direction include \cite{KowalksiSahaTsimerman2012,KnightlyReno2019,SugiyamaSuriajaya2022,Fazzari2024}.

In this article we study weighted low-lying zeros of spinor and standard $L$-functions attached to degree 2 Siegel modular forms. For a general introduction on Siegel modular forms, we refer readers to \cite{IntroductoryLecturesonSiegelModularForms,Pitale}

We proceed to describe our results. Let $k \geq 6$ be an even integer. Let $S_k(\Gamma_2)$ be the space of degree 2 holomorphic Siegel cusp forms of weight $k$ for the symplectic group $\Gamma_2 = \text{Sp}_4(\bZ)$. Each form $F \in S_k(\Gamma_2)$ is a holomorphic function on the Siegel upper half-plane
\begin{equation}
    \bH_2 = \{Z = X + iY \in M_{2}(\bC): Z = Z^T, Y > 0\},
\end{equation}
which satisfies the automorphic condition
\begin{equation}
    F((AZ+B)(CZ+D)^{-1}) = \det(CZ+D)^k F(Z), \hspace{3mm} \begin{pmatrix} A & B \\ C & D \end{pmatrix} \in \Gamma_2, \hspace{3mm} Z \in \bH_2.
\end{equation}
Here and after we use $M_n(R)$ to denote the ring of $n \times n$ matrices over a ring $R$. 

The Fourier expansion of $F$ is
\begin{equation}
    F(Z) = \sum_{T \in \cT} a_F(T) (\det T)^{\frac{k}{2}-\frac{3}{4}}e(\Tr(TZ)), \hspace{3mm} Z \in \bH_2,
\end{equation}
where the summation is taken over the set
\begin{equation}
    \cT = \{T = (t_{ij}) \in M_2(\bR): T > 0, t_{11} \in \bZ, t_{22} \in \bZ, 2t_{12} = 2t_{21} \in \bZ \}.
\end{equation}
We call $a_F(T)$ the (normalized) Fourier coefficient of $F$ at $T$. It is known that $a_F(T) \in \bR$.

We use $I$ to denote the $2 \times 2$ identity matrix. For $F \in S_k(\Gamma_2)$ we set
\begin{equation}\label{harmonic weight}
    \omega_F = \frac{\sqrt{\pi}}{4}(4\pi)^{3-2k}\Gamma(k-3/2)\Gamma(k-2)\frac{a_F(I)^2}{||F||^2}
\end{equation}
to be the "harmonic" weight attached to $F$, where $||F||$ is the Petersson norm of $F$ defined by
\begin{equation}
    ||F|| = \left(\int_{\Gamma_2 \backslash \bH_2} |F(Z)|^2 (\det Y)^k \frac{dXdY}{(\det Y)^3}\right)^\frac{1}{2}.
\end{equation}

We now choose a basis $H_k(\Gamma_2)$ of $S_k(\Gamma_2)$ consisting of eigenforms for all Hecke operators (we call such a form a Hecke eigenform). It is known (see, for example, equation (1.8) in \cite{Blomer2019}) that 
\begin{equation}\label{0th moment of harmonic weight}
    \sum_{F \in H_k(\Gamma_2)}\omega_F = 1 + O(e^{-k}).
\end{equation}

To each form $F \in H_k(\Gamma_2)$ one can attach a degree 4 spinor $L$-function $L(s,F;\text{spin})$ and a degree 5 standard $L$-function $L(s,F;\text{std})$, both normalized so that the central point is $s=1/2$. The analytic conductors of those $L$-functions are of size $k^2$ and $k^4$, respectively. Further properties of these $L$-functions are discussed in section \ref{section 2}. 

We assume GRH for both spinor and standard $L$-functions, and denote their non-trivial zeros on the critical line by
\begin{equation}
    \rho_{F,\text{spin}} = \frac{1}{2}+i\gamma_{F,\text{spin}}, \hspace{3mm} \rho_{F,\text{std}} = \frac{1}{2}+i\gamma_{F,\text{std}}.
\end{equation}
The corresponding density functions with respect to a test function $\Phi$ are
\begin{align}
    D(F;\Phi;\text{spin}) &= \sum_{\rho_{F,\text{spin}}}\Phi\left(\frac{\gamma_{F,\text{spin}}}{2\pi}\log c_{F;\text{spin}}\right), \label{density spinor L-function} \\
    D(F;\Phi;\text{std}) &= \sum_{\rho_{F,\text{std}}}\Phi\left(\frac{\gamma_{F,\text{std}}}{2\pi}\log c_{F;\text{std}}\right). \label{density standard L-function}
\end{align}

Our first result concerning low-lying zeros of spinor $L$-functions is as follows:
\begin{theorem}\label{spin L-function}
    Let $\Phi$ be an even Schwartz function whose Fourier transform has support in $(-1,1)$. For $F \in H_k(\Gamma_2)$, define $D(F;\Phi;\text{spin})$ as in \ref{density spinor L-function} with $c_{F;\text{spin}} = k^2$ and $\omega_F$ as in \ref{harmonic weight}. Assume GRH for $L(s,F;\text{spin})$. Then we have
    \begin{equation}
        \lim_{k \to \infty} \sum_{F \in H_k(\Gamma_2)} \omega_F D(F;\Phi;\text{spin}) = \hat{\Phi}(0) - \frac{\Phi(0)}{2} = \int_{-\infty}^\infty \Phi(x)W(Sp)(x)dx.
    \end{equation}
\end{theorem}

\begin{remark}
    The result above has been obtained in \cite{KowalksiSahaTsimerman2012}, but only for test functions $\Phi$ with supp$(\hat{\Phi}) \subset (-\frac{4}{15},\frac{4}{15})$, as an application of their quantitative local equidistribution result. Here we extend the range of support to $(-1,1)$. This improvement is crucial in application to non-vanishing problems, as we will explain in section \ref{section 5}.
\end{remark}

Let $H_k^*(\Gamma_2) \subset H_k(\Gamma_2)$ denote a Hecke basis of the space of Saito-Kurokawa lifts (these concepts will be discussed in section \ref{section 2}). Then as a direct corollary of Theorem \ref{spin L-function} we can establish the following non-vanishing result:
\begin{corollary}\label{non-vanishing}
    Assume GRH for $L(s,F;\text{spin})$. Then we have
    \begin{equation}
        \liminf_{k \to \infty} \sum_{\substack{F \in H_k(\Gamma_2) \backslash H_k^*(\Gamma_2) \\ L(1/2,F;\text{spin}) \neq 0}}\omega_F \geq \frac{3}{4}.
    \end{equation}
\end{corollary}

\begin{remark}
    For comparison, it is shown in \cite{Blomer2019} that 
    \begin{equation}
        \sum_{\substack{F \in H_k(\Gamma_2) \backslash H_k^*(\Gamma_2) \\ L(1/2,F;\text{spin}) \neq 0}}\omega_F \gg \frac{1}{\log k},
    \end{equation}
    unconditionally for large $k$. This follows from asymptotic formulas for the first and second moments of central values. Although it is not surprising that GRH would yield much stronger result, one still needs the range of support in Theorem \ref{spin L-function} not to be too small to carry out the argument.
\end{remark}

For low-lying zeros of standard $L$-functions, we have the following result:
\begin{theorem}\label{standard L-function}
   Let $\Phi$ be an even Schwartz function whose Fourier transform has support in $(-\frac{1}{4},\frac{1}{4})$. For $F \in H_k(\Gamma_2)$, define $D(F;\Phi;\text{std})$ as in \ref{density standard L-function} with $c_{F;\text{std}} = k^4$ and $\omega_F$ as in \ref{harmonic weight}. Assume GRH for $L(s,F;\text{std})$. Then we have
    \begin{equation}
        \lim_{k \to \infty} \sum_{F \in H_k(\Gamma_2)} \omega_F D(F;\Phi;\text{std}) = \hat{\Phi}(0) - \frac{\Phi(0)}{2} = \int_{-\infty}^\infty \Phi(x)W(Sp)(x)dx.
    \end{equation}
\end{theorem}

\begin{remark}
    An unweighted version of Theorem \ref{standard L-function} was established in \cite{KimWakatsukiYamauchi2020}, for test functions $\Phi$ whose Fourier transform have sufficiently small support (for a precise range of support, see Proposition 9.3 in \cite{KimWakatsukiYamauchi2020}). The (unweighted) symmetry type is also symplectic. For comparison, the symmetry type of low-lying zeros of spinor $L$-functions changes from orthogonal to symplectic when weighted by $\omega_F$. 
\end{remark}

We may further extend the range of support in Theorem \ref{standard L-function} from $(-\frac{1}{4},\frac{1}{4})$ to $(-\frac{5}{18},\frac{5}{18})$ by performing an extra (smooth) average over weight $k$. Our result is the following:
\begin{theorem}\label{average standard L-function}
     Let $\Omega \in C_c^\infty(0,\infty)$ be such that $\Omega \geq 0$, not identically 0. Let $\Phi$ be an even Schwartz function whose Fourier transform has support in $(-\frac{5}{18},\frac{5}{18})$. For $F \in H_k(\Gamma_2)$ and large parameter $K > 0$, define $D(F;\Phi;\text{std})$ as in \ref{density standard L-function} with $c_{F;\text{std}} = K^4$ and $\omega_F$ as in \ref{harmonic weight}. Assume GRH for $L(s,F;\text{std})$. Then we have
    \begin{equation}
        \lim_{K \to \infty}\left(\sum_k \Omega\left(\frac{k}{K}\right)\right)^{-1} \sum_k \Omega \left(\frac{k}{K}\right) \sum_{F \in H_k(\Gamma_2)} \omega_F D(F;\Phi;\text{std}) = \hat{\Phi}(0) - \frac{\Phi(0)}{2} = \int_{-\infty}^\infty \Phi(x)W(Sp)(x)dx.
    \end{equation}
    where the summation in $k$ is over even integers. 
\end{theorem}

This article is organized as follows: In section \ref{section 2}, we first review some facts about spinor and standard $L$-functions. We then work out the combinatorial relations between certain functions in Satake parameters of a form $F \in H_k(\Gamma_2)$ and its Fourier coefficients at scalar matrices. These relations allow us to apply the Petersson formula, which we state in section \ref{section 3}. In section \ref{section 3} we also take average over weight $k$ in the Petersson formula and give an upper bound for the off-diagonal term. In section \ref{section 4} we apply the results established in previous sections, as well as the explicit formula to prove Theorems \ref{spin L-function}-\ref{average standard L-function}. In section \ref{section 5} we prove Corollary \ref{non-vanishing} and discuss some other issues concerning non-vanishing of central $L$-values. 

\section{Spinor and Standard L-functions}\label{section 2}
Let $F \in H_k(\Gamma_2)$ be a Hecke eigenform. It is known that for each prime $p$ there are 3 complex numbers $\alpha_{F,0}(p),\alpha_{F,1}(p),\alpha_{F,2}(p)$, called the Satake parameters of $F$ at $p$, with certain prescribed properties. See chapter 3 of \cite{Pitale} for a detailed discussion. In particular, these Satake parameters satisfy the relation
\begin{equation}\label{Satake relation}
    \alpha_{F,0}(p)^2\alpha_{F,1}(p)\alpha_{F,2}(p) = 1.
\end{equation}

Let $S_{2k-2}(\Gamma_1)$ denote the space of holomorphic cusp forms of weight $2k-2$ for the full modular group $\Gamma_1 =$ SL$_2(\bZ)$. There is an injective Hecke-equivariant linear map 
\begin{equation}
    SK: S_{2k-2}(\Gamma_1) \to S_k(\Gamma_2), \hspace{3mm} f \mapsto F_f,
\end{equation}
called the Saito-Kurokawa lifting. We denote the image of $SK$ by $S_k^*(\Gamma_2)$ and call forms in $S_k^*(\Gamma_2)$ Saito-Kurokawa lifts. We also use $H_k^*(\Gamma_2)$ for a basis of $S_k^*(\Gamma_2)$ consisting of Hecke eigenforms. There are various ways to construct such a lifting map. For a construction using half-integral weight modular forms, see section 2.1.3 of \cite{Pitale} and references there. 

The Generalized Ramanujan Conjecture (GRC) asserts that all Satake parameters $\alpha_{F,i}(p)$ ($i=0,1,2$) have absolute value 1. For elliptic modular forms, the Ramanujan Conjecture is true due to Deligne's work \cite{Deligne1974}. For degree 2 Siegel modular forms $F \in H_k(\Gamma_2)$, GRC is proved to be true only for non-Saito-Kurokawa lifts, due to R. Weissauer's result \cite{Weissauer2009}. For Saito-Kurokawa lifts $F_f \in H_k^*(\Gamma_2)$, GRC is false, as we will see in Andrianov's explicit formula \ref{Andrianov explicit formula} below. 

\subsection{The Spinor L-function}
The spinor $L$-function attached to a Hecke eigenform $F \in H_k(\Gamma_2)$ is defined by a degree 4 Euler product
\begin{equation}
\begin{split}
    L(s,F;\text{spin}) &= \prod_p\left(1-\frac{\alpha_{F,0}(p)}{p^s}\right)^{-1}\left(1-\frac{\alpha_{F,0}(p)\alpha_{F,1}(p)}{p^s}\right)^{-1}\left(1-\frac{\alpha_{F,0}(p)\alpha_{F,2}(p)}{p^s}\right)^{-1} \\
    &\times \left(1-\frac{\alpha_{F,0}(p)\alpha_{F,1}(p)\alpha_{F,2}(p)}{p^s}\right)^{-1},
\end{split}
\end{equation}
which converges absolutely in some right half-plane. By setting 
\begin{equation}
\alpha_F(p) = \alpha_{F,0}(p), \hspace{3mm} \beta_F(p) = \alpha_{F,0}(p)\alpha_{F,1}(p),
\end{equation}
we may rewrite the above Euler product as
\begin{equation}\label{spinor L-function Euler product}
    L(s,F;\text{spin}) = \prod_p\left(1-\frac{\alpha_{F}(p)}{p^s}\right)^{-1}\left(1-\frac{\beta_F(p)}{p^s}\right)^{-1}\left(1-\frac{\alpha_F(p)^{-1}}{p^s}\right)^{-1}\left(1-\frac{\beta_F(p)^{-1}}{p^s}\right)^{-1},
\end{equation}
in view of the relation \ref{Satake relation}. 

It is proved by A. N. Andrianov \cite{Andrianov1974} that $L(s,F;\text{spin})$ extends to a meromorphic function on $\bC$, which has a simple pole at $s=\frac{3}{2}$ if $F$ is a Saito-Kurokawa lift, and is entire otherwise. Its functional equation takes the following form:
\begin{equation}\label{spin L-function functional equation}
    \Lambda(s,F;\text{spin}) = \Gamma_\bC(s+1/2)\Gamma_\bC(s+k-3/2)L(s,F;\text{spin}) = \Lambda(1-s,F;\text{spin}),
\end{equation}
where $\Gamma_\bC(s) = 2(2\pi)^{-s}\Gamma(s)$. For $F_f \in H_k^*(\Gamma_2)$ a Saito-Kurokawa lift, its spinor $L$-function decomposes as follows:
\begin{equation}
    L(s,F_f;\text{spin}) = \zeta(s+1/2)\zeta(s-1/2)L(s,f),
\end{equation}
where $L(s,f)$ is the Hecke $L$-function of the elliptic cusp form $f$. 

For $F \in H_k(\Gamma_2)$ we have the following Andrianov's explict formula \cite{Andrianov1974}:
\begin{equation}\label{Andrianov explicit formula}
    a_F(I)L(s,F;\text{spin}) = \zeta(s+1/2)L(s+1/2,\chi_{-4})\sum_{n=1}^\infty \frac{a_F(nI)}{n^s},
\end{equation}
where $\chi_{-4}$ is the non-trivial Dirichlet character modulo 4. From this formula it follows
\begin{equation}\label{observation}
    a_F(I) = 0 \implies a_F(nI) = 0, \hspace{3mm} n \geq 1.
\end{equation}

\subsection{The Standard L-function}
The standard $L$-function attached to a Hecke eigenform $F \in H_k(\Gamma_2)$ is defined by a degree 5 Euler product 
\begin{equation}
L(s,F;\text{std}) = \prod_p \left(1-\frac{1}{p^s}\right)^{-1}\left(1-\frac{\alpha_{F,1}(p)}{p^s}\right)^{-1}\left(1-\frac{\alpha_{F,1}(p)^{-1}}{p^s}\right)^{-1}\left(1-\frac{\alpha_{F,2}(p)}{p^s}\right)^{-1}\left(1-\frac{\alpha_{F,2}(p)^{-1}}{p^s}\right)^{-1},
\end{equation}
which converges absolutely in some right half-plane. Using \ref{Satake relation}, we rewrite this Euler product as
\begin{equation}
\begin{split}
    L(s,F;\text{std}) &= \prod_p \left(1-\frac{1}{p^s}\right)^{-1}\left(1-\frac{\alpha_F(p)\beta_F(p)}{p^s}\right)^{-1}\left(1-\frac{\alpha_F(p)^{-1}\beta_F(p)}{p^s}\right)^{-1} \\
    &\times \left(1-\frac{\alpha_F(p)\beta_F(p)^{-1}}{p^s}\right)^{-1}\left(1-\frac{\alpha_F(p)^{-1}\beta_F(p)^{-1}}{p^s}\right)^{-1}.
\end{split}
\end{equation}

The analytic continuation and functional equation of standard $L$-functions were worked out by S. B\"ocherer \cite{Bocherer1985}. He proved that $L(s,F;\text{std})$ extends to an entire function and satisfies a functional equation
\begin{equation}
    \Lambda(s,F;\text{std}) = \Gamma_\bR(s)\Gamma_\bC(s+k-1)\Gamma_\bC(s+k-2)L(s,F;\text{std}) = \Lambda(1-s,F;\text{std}),
\end{equation}
where $\Gamma_\bR = \pi^{-\frac{s}{2}}\Gamma(s/2)$ and $\Gamma_\bC(s) = 2(2\pi)^{-s}\Gamma(s)$.

\subsection{Combinatorial Relations}
For $F \in H_k(\Gamma_2)$, $m \geq 1$ and prime $p$, we set
\begin{align}
    c_m(p;F) &= \alpha_F(p)^m + \alpha_F(p)^{-m} + \beta_F(p)^m + \beta_F(p)^{-m}, \label{power sum spinor} \\
    \tau_{2m}(p;F) &= 1 + \alpha_F(p)^m \beta_F(p)^m + \alpha_F(p)^m \beta_F(p)^{-m} + \alpha_F(p)^{-m} \beta_F(p)^m + \alpha_F(p)^{-m} \beta_F(p)^{-m} \label{power sum standard}
\end{align}
to be the $m^{th}$-power sum of local parameters of $L(s,F;\text{spin})$ and $L(s,F;\text{std})$ at $p$ respectively.

The main goal of this section is to find expressions of these power sums in terms of Fourier coefficients of $F$ at scalar matrices, for $m = 1,2$, under the assumption that $a_F(I) \neq 0$. Note that the condition $a_F(I) \neq 0$ is not a direct consequence of $F \neq 0$, unlike in the elliptic case, where a primitive form $f$ vanishes if and only if its first Fourier coefficient vanishes. In fact, determining whether $a_F(I) = 0$ or not is a difficult problem because $\omega_F$ is intimately connected to central values of spinor $L$-functions (B\"ocherer's conjecture). However, as we shall see later in section \ref{section 4}, making this assumption here does no harm to our argument. Our result is as follows: 

\begin{lemma}\label{combinatorial lemma}
    Let $F \in H_k(\Gamma_2)$ be a Hecke eigenform. For any prime $p$ and $m \geq 1$, define $c_m(p;F)$ and $\tau_{2m}(p;F)$ as in \ref{power sum spinor} and \ref{power sum standard}. Assume that $a_F(I) \neq 0$, and set $U_m(p;F) = \frac{a_F(p^mI)}{a_F(I)}$. Also set $\lambda_p = 1 + \chi_{-4}(p)$ and $\mu_p = \chi_{-4}(p)$, where $\chi_{-4}$ is the non-trivial Dirichlet character modulo 4. Then we have
    \begin{align*}
        c_1(p;F) &= U_1(p;F) + \frac{\lambda_p}{\sqrt{p}}, \\
        c_2(p;F) &= - U_1(p;F)^2 + 2U_2(p;F) + \frac{\lambda_p^2 - 2\mu_p}{p}, \\
        \tau_2(p;F) &= U_1(p;F)^2 -U_2(p;F)+ \frac{\lambda_p}{\sqrt{p}}U_1(p;F) + \frac{\mu_p}{p} - 1, \\
        \tau_4(p;F) &= -U_3(p;F)U_1(p;F) + U_2(p;F)^2 + \frac{\lambda_p}{\sqrt{p}}U_2(p;F)U_1(p;F) + \left(\frac{\mu_p}{p} - 1\right)U_1(p;F)^2 \\
        &- \frac{\lambda_p}{\sqrt{p}}U_3(p;F) + \left(\frac{\lambda_p^2 - 2\mu_p}{p}\right)U_2(p;F) + \left(\frac{\lambda_p\mu_p}{p^\frac{3}{2}} - 2\frac{\lambda_p}{\sqrt{p}}\right)U_1(p;F) + \frac{\mu_p^2}{p^2} - \frac{\lambda_p^2}{p} + 1.
    \end{align*}
\end{lemma}

\begin{remark}
    The key feature of Lemma \ref{combinatorial lemma} is that we are able to express $\tau_4(p;F)$ using polynomials of $U_m(p;F)$ of degree 2 (and not of higher degree). This is essential when we deal with weighted low-lying zeros of standard $L$-functions using the Petersson formula. 
\end{remark}

\begin{proof}
    Throughout the proof, $F$ and $p$ are fixed. To save notations we use $c_m,\tau_{2m},U_m,\alpha_p,\beta_p$ to denote $c_m(p;F),\tau_{2m}(p;F),U_m(p;F),\alpha_F(p),\beta_F(p)$ respectively, with the understanding that they depend on $F$ and $p$. 
    
    We start with Andrianov's explicit formula \ref{Andrianov explicit formula}:
    $$a_F(I)L(s,F;\text{spin}) = \zeta(s+1/2)L(s+1/2,\chi_{-4})\sum_{n=1}^\infty \frac{a_F(nI)}{n^s}.$$
    Using Euler product expansions for the $L$-functions involved, we see that the two Dirichlet series 
    \begin{equation}
        \prod_p\left(1-\frac{\alpha_p}{p^s}\right)\left(1-\frac{\beta_p}{p^s}\right)\left(1-\frac{\alpha_p^{-1}}{p^s}\right)\left(1-\frac{\beta_p^{-1}}{p^s}\right)\sum_{n=1}^\infty \frac{a_F(nI)a_F(I)^{-1}}{n^s}
    \end{equation}
    and 
    \begin{equation}
        \prod_p\left(1-\frac{1}{p^{s+1/2}}\right)\left(1-\frac{\chi_{-4}(p)}{p^{s+1/2}}\right)
    \end{equation}
    both converge absolutely in some right half-plane and are equal. Comparing coefficients of $p^{-as}$, $a = 1,2,3,4$, we obtain
    \begin{align}
        -\frac{\lambda_p}{\sqrt{p}} &= U_1 - c_1, \label{relation 1}  \\
        \frac{\mu_p}{p} &= U_2 - U_1c_1 + \tau_2 + 1, \label{relation 2} \\
        0 &= U_3 - U_2c_1 + U_1(\tau_2 + 1) - c_1,  \label{relation 3} \\
        0 &= U_4 - U_3c_1 + U_2(\tau_2 + 1) - U_1c_1 + 1. \label{relation 4}
    \end{align}
    
    We also have elementary relations
    \begin{align}
        c_1^2 &= c_2 + 2(\tau_2 + 1), \label{elementary relation 1} \\
        (\tau_2+1)^2 &= 3 + \tau_4 + 4\tau_2 + 2c_2. \label{elementary relation 2}
    \end{align}
    
    From equations \ref{relation 1} and \ref{relation 2} we obtain directly
    \begin{align}
        c_1 &= U_1 + \frac{\lambda_p}{\sqrt{p}}, \label{c_1} \\
        \tau_2 &= U_1^2 - U_2 + \frac{\lambda_p}{\sqrt{p}}U_1 + \frac{\mu_p}{p} - 1. \label{tau_2}
    \end{align}
    
    Combining equations \ref{elementary relation 1}, \ref{c_1} and \ref{tau_2} we have
    \begin{equation}\label{c_2}
        c_2 = -U_1^2 + 2U_2 + \frac{\lambda_p^2 - 2\mu_p}{p}.
    \end{equation}
    
    Using equations \ref{elementary relation 2}, \ref{tau_2} and \ref{c_2} we express $\tau_4$ as
    \begin{equation}\label{pre tau_4}
    \begin{split}
    \tau_4 &= U_1^4 - 2U_2U_1^2 + 2\frac{\lambda_p}{\sqrt{p}}U_1^3 + U_2^2 -2\frac{\lambda_p}{\sqrt{p}}U_2U_1 + \left(\frac{\lambda_p^2+2\mu_p}{p}-2\right)U_1^2 \\
     &- 2\frac{\mu_p}{p}U_2 + \left(2\frac{\lambda_p\mu_p}{p^\frac{3}{2}} - 4\frac{\lambda_p}{\sqrt{p}}\right)U_1 + \frac{\mu_p^2}{p^2} - 2\frac{\lambda_p^2}{p} + 1.
    \end{split}
    \end{equation}
    However, to get the final form of $\tau_4$, we must express $U_1^4$, $U_2U_1^2$ and $U_1^3$ using degree 2 polynomials in $U_a (a = 1,2,3,4).$ Combining equations \ref{relation 3}, \ref{c_1} and \ref{tau_2}, we have
    \begin{equation}\label{U_1^3}
        U_1^3 = 2U_2U_1 - \frac{\lambda_p}{\sqrt{p}}U_1^2 - U_3 + \frac{\lambda_p}{\sqrt{p}}U_2 + \left(1-\frac{\mu_p}{p}\right)U_1 + \frac{\lambda_p}{\sqrt{p}}. 
    \end{equation}
    Likely, equations \ref{relation 4}, \ref{c_1} and \ref{tau_2} give us
    \begin{equation}\label{U_2U_1^2}
        U_2U_1^2 = U_3U_1 + U_2^2 - \frac{\lambda_p}{\sqrt{p}}U_2U_1 + U_1^2 - U_4 + \frac{\lambda_p}{\sqrt{p}}U_3 - \frac{\mu_p}{p}U_2 + \frac{\lambda_p}{\sqrt{p}}U_1 - 1.
    \end{equation}
    Further, we multiply \ref{U_1^3} by $U_1$ and apply \ref{U_1^3}, \ref{U_2U_1^2} to get
    \begin{equation}\label{U_1^4}
    \begin{split}
        U_1^4 &= U_3U_1 + 2U_2^2 - 3\frac{\lambda_p}{\sqrt{p}}U_2U_1 + \left(3+\frac{\lambda_p^2-\mu_p}{p}\right)U_1^2 - 2U_4 + 3\frac{\lambda_p}{\sqrt{p}}U_3 \\
         &- \left(\frac{\lambda_p^2 + 2\mu_p}{p}\right)U_2 + \left(2\frac{\lambda_p}{\sqrt{p}} + \frac{\lambda_p\mu_p}{p^\frac{3}{2}}\right)U_1 - \left(2 + \frac{\lambda_p^2}{p}\right).
    \end{split}
    \end{equation}
    Finally, we insert \ref{U_1^3}, \ref{U_2U_1^2} and \ref{U_1^4} into \ref{pre tau_4} to get
    \begin{equation}\label{tau_4}
    \begin{split}
        \tau_4 &= -U_3U_1 + U_2^2 + \frac{\lambda_p}{\sqrt{p}}U_2U_1 + \left(\frac{\mu_p}{p} - 1\right)U_1^2 - \frac{\lambda_p}{\sqrt{p}}U_3 + \left(\frac{\lambda_p^2 - 2\mu_p}{p}\right)U_2 \\
        &+ \left(\frac{\lambda_p\mu_p}{p^\frac{3}{2}} - 2\frac{\lambda_p}{\sqrt{p}}\right)U_1 + \frac{\mu_p^2}{p^2} - \frac{\lambda_p^2}{p} + 1.
    \end{split}
    \end{equation}
\end{proof}

\section{The Petersson Formula}\label{section 3}
The main tool used in this paper is a spectral summation formula of Petersson type. This formula was first worked by Y. Kitaoka \cite{Kitaoka1984} by computing Fourier coefficients of Siegel Poincar\'e series. In this section we introduce this formula and consider an averaged (over weight) version of it. 

We begin by introducing some notations. For $k \geq 6$ even, we set 
\begin{equation}
    c_k = \frac{\sqrt{\pi}}{4}(4\pi)^{3-2k}\Gamma(k-3/2)\Gamma(k-2).
\end{equation}

For $T,Q \in \cT$ we define 
\begin{equation}\label{spectral side}
    \Delta_k(T,Q) = \sum_{F \in H_k(\Gamma_2)} c_k \frac{a_F(T)a_F(Q)}{||F||^2}.
\end{equation}

For a matrix $C \in M_2(\bZ)$ with $\det C \neq 0$ (we denote the set of such matrices by $\cC$) and $Q,T \in \cT$, define the symplectic Kloosterman sum to be
\begin{equation}
    K(Q,T;C) = \sum_D e(\Tr(AC^{-1}Q + C^{-1}DT)),
\end{equation}
where $D$ runs through the set 
\begin{equation}
    \{D \in M_2(\bZ) \bmod C\Lambda: \begin{pmatrix}
        A & * \\ C & D
    \end{pmatrix} \in \Gamma_2 \},
\end{equation}
and $\Lambda$ is the set of $2 \times 2$ symmetric integral matrices. By elementary divisor theory one has the following estimate \cite{Kitaoka1984}:
\begin{equation}\label{symplectic Kloosterman sum trivial bound}
    |K(Q,T;C)| \leq |\det C|^{3/2}.
\end{equation}

For $P = \begin{pmatrix} p_1 & p_2/2 \\ p_2/2 & p_4 \end{pmatrix}, S = \begin{pmatrix} s_1 & s_2/2 \\ s_2/2 & s_4 \end{pmatrix} \in \cT$ and $c \geq 1$, we define another exponential sum
\begin{equation}
    H^{\pm}(P,S;c) = \delta_{s_4=p_4}\sideset{}{^*}\sum_{d_1 \bmod c}\sum_{d_2 \bmod c}e\left(\frac{\overline{d_1}s_4d_2^2 \mp\overline{d_1}p_2d_2+s_2d_2+\overline{d_1}p_1 + d_1s_1}{c}\mp\frac{p_2s_2}{2cs_4}\right).
\end{equation}
For these we have the trivial bound
\begin{equation}\label{trivial estimate H}
    |H^{\pm}(P,S;c)| \leq c^2.
\end{equation}

For $P \in M_2(\bR)$ with positive eigenvalues $\lambda_1,\lambda_2 > 0$ we set 
\begin{equation}
    \cJ_{k-3/2}(P) =\int_0^\frac{\pi}{2} J_{k-3/2}(4\pi \sqrt{\lambda_1}\sin \theta)J_{k-3/2}(4\pi \sqrt{\lambda_2}\sin \theta)\sin \theta d\theta,
\end{equation}
where $J_{k-3/2}$ is the usual $J$-Bessel function of half-integral order $k-\frac{3}{2}$. With these notations, we can now state the Petersson formula.
\begin{lemma}\label{Petersson formula lemma}
For $T,Q \in \cT$ and $k \geq 6$ even, define $\Delta_k(T,Q)$ as in \ref{spectral side}. Then we have
    \begin{equation}\label{Petersson formula}
        \Delta_k(T,Q) = \frac{1}{8}|\text{Aut}(T)|\left(\frac{\det Q}{\det T}\right)^{\frac{k}{2}-\frac{3}{4}}\delta_{Q \sim T} + \frac{\sqrt{2}\pi}{8}G_{1,k}
        (T,Q) + \pi^2G_{2,k}(T,Q),
    \end{equation}
    where $\text{Aut}(T), G_{1,k}(T,Q)$ and $G_{2,k}(T,Q)$ are defined by
    \begin{align}
        \text{Aut}(T) &= \{U \in \text{GL}_2(\bZ): U^TTU = T\}, \\
        G_{1.k}(T,Q) &= \sum_{\pm}\sum_{s=1}^\infty\sum_{c=1}^\infty\sum_{U,V}\frac{(-1)^{\frac{k}{2}}}{c^\frac{3}{2}s^\frac{1}{2}}H^{\pm}(UQU^T,V^{-1}TV^{-T};c)J_{k-3/2}\left(\frac{4\pi\sqrt{\det(TQ)}}{cs}\right), \label{rank 1 term}\\
        G_{2,k}(T,Q) &= \sum_{C \in \cC}\frac{K(Q,T;C)}{|\det C|^\frac{3}{2}} \cJ_{k-3/2}(TC^{-1}QC^{-T}). \label{rank 2 term}
    \end{align}
    Here the summation $\sum_{U,V}$ in \ref{rank 1 term} is over $U = (u_{ij}) / \{\pm I\},V = (v_{ij}) \in $ GL$_2(\bZ)$ such that 
    \begin{equation}
        (u_{21},u_{22})Q(u_{21},u_{22})^T = (-v_{21},v_{11})T(-v_{21},v_{11})^T = s.
    \end{equation}
    The delta symbol $\delta_{Q \sim T}$ is equal to 1 if $Q$ and $T$ are equivalent in the sense of quadratic forms, and is equal to 0 otherwise.  
\end{lemma}

\begin{remark}
    Following Kitaoka \cite{Kitaoka1984}, we call the three terms in \ref{Petersson formula} containing $\delta_{Q \sim T}$, $G_{1,k}(T,Q)$ and $G_{2,k}(T,Q)$ the diagonal term, the rank 1 term and the rank 2 term respectively. Note that the classical Petersson formula for elliptic modular forms contains only a diagonal term and an off-diagonal term. 
\end{remark}

\begin{remark}
    As pointed out by V. Blomer (see Remark 1 in \cite{Blomer2019}), there are some numerical errors in Kitaoka's original derivation of the Petersson formula. The version that we present here is based on Lemma 1 in \cite{Blomer2019}. However, our results do not depend on exact values of those constants. 
\end{remark}

The main purpose of this section is to establish the following averaged Petersson formula:
\begin{lemma}\label{averaged Petersson}
    Let $m,n \geq 1$ be positive integers such that $m|n$. For $k \geq 6$ even, define $\Delta_{k}(mI,nI)$ as in \ref{spectral side}. Let $\Omega \in C_c^\infty(0,\infty)$ be such that $\Omega \geq 0$, not identically 0. Then for large $K > 0$ we have 
    \begin{equation}\label{averaged Petersson formula}
        \left(\sum_k \Omega\left(\frac{k}{K}\right)\right)^{-1} \sum_k \Omega\left(\frac{k}{K}\right)\Delta_k(mI,nI) = \delta_{m=n} + O_{j,\epsilon,\Omega}\left(\left(\frac{m^{\frac{3}{2}-\epsilon}n^{-\frac{1}{2}+\epsilon}}{K^4}\right) + \left(\frac{(mn)^{2+\epsilon}}{K^{5+2\epsilon}}\right) + \left(\frac{(mn)^{\frac{j}{2}+1}}{K^{2j+3}}\right)\right)
    \end{equation}
for any $j \geq 3$ and $\epsilon > 0$ small. Here the summation $\sum_{k}$ is over positive even integers $k \geq 6$. 
\end{lemma}

\begin{remark}
    Our Lemma \ref{averaged Petersson} can be viewed as a GSp$_4$ analogue of the classical averaged Petersson formula on GL$_2$. See equation (5.81) in \cite{TopicsinClassicalAutomorphicForms}. The main difficulty is the presence of a product of two Bessel functions (instead of a single Bessel function), each of half-integral order (instead of integral order). As we shall see in the proof below, this can be overcome by applying an integral representation \ref{product formula} of a product of two Bessel functions. 
\end{remark}

\begin{proof}
    After applying the Petersson formula \ref{Petersson formula}, we divide the left side of \ref{averaged Petersson formula} into three terms. We also set $g(x) = \Omega(\frac{x}{K})$ and $\ell = k-3/2$ to save notations. 

    We denote the contribution of the diagonal term by $R_0$. Thus
    \begin{equation}
        R_0 = \frac{1}{8}\left(\sum_{k}g(k)\right)^{-1}\sum_k g(k) |\text{Aut}(mI)|\left(\frac{n}{m}\right)^{\ell}\delta_{mI \sim nI}.
    \end{equation}
    Note that $mI$ and $nI$ define the same quadratic form if and only if $m=n$, and that $|\text{Aut}(mI)| = 8$. Thus the above expression reduces to $R_0 = \delta_{m=n}$.

    Denote by $R_1$ the sum of the rank 1 term over $k$. We have 
    \begin{equation}\label{pre R_1}
        R_1 = \sum_k g(k) \sum_{\pm}\sum_{s=1}^\infty\sum_{c=1}^\infty\sum_{U,V}\frac{(-1)^{\frac{k}{2}}}{c^\frac{3}{2}s^\frac{1}{2}}H^{\pm}(nUU^T,mV^{-1}V^{-T};c)J_{\ell}\left(\frac{4\pi mn}{cs}\right),
    \end{equation}
    where the sum $\sum_{U,V}$ is over 
    \begin{equation}
        n(u_{21}^2 + u_{22}^2) = m(v_{11}^2+v_{21}^2) = s.
    \end{equation}
    So in particular $n|s$. Making change of variable $s \mapsto ns$, we may rewrite $R_1$ as 
    \begin{equation}\label{R_1}
        R_1 = \sum_k g(k) \sum_{\pm}\sum_{s=1}^\infty\sum_{c=1}^\infty\sum_{U,V}\frac{(-1)^{\frac{k}{2}}}{c^\frac{3}{2}(ns)^\frac{1}{2}}H^{\pm}(nUU^T,mV^{-1}V^{-T};c)J_{\ell}\left(\frac{4\pi m}{cs}\right),
    \end{equation}
    where $\sum_{U,V}$ is over 
    \begin{equation}\label{number of solutions}
        u_{21}^2 + u_{22}^2 = s, \hspace{3mm} v_{11}^2 + v_{22}^2 = \frac{n}{m}s.
    \end{equation}
    These equations have $O(s^\epsilon)$ and $O((\frac{n}{m}s)^\epsilon)$ integral solutions, respectively, for any $\epsilon>0$, by the fact
    \begin{equation}\label{number of integral solutions}
        |\{(x,y) \in \bZ^2: x^2 + y^2 = s\}| = O(s^\epsilon).
    \end{equation}
    
    In view of the estimate \cite{GradshteynRyzhik2015}
    \begin{equation}\label{Bessel estimate 1}
        J_{\ell}(x) \ll \left(\frac{x}{\ell}\right)^\ell, \hspace{3mm} x > 0, \ell > \frac{1}{2},
    \end{equation}
    we may cut-off the sum in \ref{R_1} by $sc \ll \frac{m}{K}$ up to a negligible error. In this range we change summation order and deal with the inner sum
    \begin{equation}
        \sum_k g(k)(-1)^{\frac{k}{2}}J_{\ell}\left(\frac{4\pi m}{cs}\right)
    \end{equation}
    by applying a lemma due to Blomer and A. Corbett (Lemma 20 in \cite{BlomerCorbett2022}) to obtain
    \begin{equation}
        \sum_kg(k)(-1)^\frac{k}{2}J_{\ell}\left(\frac{4\pi m}{cs}\right) = \omega_0\left(\frac{4\pi m}{cs}\right) + e^{\frac{4\pi i m}{cs}}\omega_+\left(\frac{4\pi m}{cs}\right) + e^{-\frac{4\pi i m}{cs}}\omega_-\left(\frac{4\pi m}{cs}\right),
    \end{equation}
    where $\omega_0(x),\omega_{\pm}(x)$ are some smooth functions on $(0,\infty)$ satisfying 
    \begin{align}
        \omega_0(x) &\ll_A K^{-A}, \\
        \omega_{\pm}(x) &\ll_A \left(1+\frac{K^2}{x}\right)^{-A},
    \end{align}
    for any $A > 0$. The contribution of the $\omega_0$ term is negligible, while the contribution of $\omega_{\pm}$ term depends on the size of $x = \frac{4\pi m}{cs}$. For example, for $x \leq K^2$ (i.e. $cs \gg \frac{m}{K^2}$), we have 
    \begin{equation}
        \omega_{+}(x) \ll_A \left(1+\frac{K^2}{x}\right)^{-A} \leq \left(\frac{K^2}{x}\right)^{-A} \ll_A  K^{-2A}m^A(cs)^{-A},
    \end{equation}
    for any $A > 0$. This estimate, together with \ref{trivial estimate H} and \ref{number of integral solutions}, give rise to
    \begin{equation}
    \begin{split}
        R_1^{+,cs \gg \frac{m}{K^2}} &=  \sum_{\pm}\sum_{\frac{m}{K^2} \ll sc \ll \frac{m}{K}}\sum_{U,V}\frac{1}{c^\frac{3}{2}(ns)^\frac{1}{2}}H^{\pm}(nUU^T,mV^{-1}V^{-T};c)\omega_+\left(\frac{4\pi m}{cs}\right) \\
        &\ll_{\epsilon,A} \sum_{\frac{m}{K^2} \ll sc\ll \frac{m}{K}} c^{-\frac{3}{2}}(ns)^{-\frac{1}{2}}s^\epsilon\left(\frac{n}{m}s\right)^\epsilon c^2 K^{-2A} m^A c^{-A}s^{-A} \\
        &\ll_{\epsilon} m^{\frac{3}{2}-\epsilon} n^{-\frac{1}{2}+\epsilon} K^{-3},
    \end{split}
    \end{equation}
    for any small $\epsilon>0$ if one fixes some $A > 3/2$. The case where $cs \ll \frac{m}{K^2}$ is analyzed similarly, and its contribution $R_1^{+,cs \ll \frac{m}{K^2}}$ to  is again at most $m^{\frac{3}{2}-\epsilon}n^{-\frac{1}{2}+\epsilon}K^{-3}$. Therefore, we have obtained that 
    \begin{equation}
       R_1 \ll_\epsilon m^{\frac{3}{2}-\epsilon}n^{-\frac{1}{2}+\epsilon}K^{-3},
    \end{equation}
    for any small $\epsilon > 0$.

    Denote by $R_2$ the sum of the rank 2 term over $k$. Explicitly, 
    \begin{equation}\label{R_2}
        R_2 = \sum_k g(k) \sum_{C \in \cC}\frac{K(nI,mI;C)}{|\det C|^\frac{3}{2}}\int_0^\frac{\pi}{2}J_{\ell}(4\pi\sqrt{\lambda_1}\sin \theta)J_{\ell}(4\pi\sqrt{\lambda_2}\sin \theta)\sin \theta d\theta,
    \end{equation}
    where $\lambda_1,\lambda_2$ are eigenvalues of the matrix $mnC^{-1}C^{-T}$. We set $\lambda_{\min}$ and $\lambda_{\max}$ to be the smaller and the larger eigenvalue of $mnC^{-1}C^{-T}$ respectively. Denote by $||\cdot||_F$ the Frobenius matrix norm. Then by Lemma 2 in \cite{Blomer2019} we have 
    \begin{equation}\label{smaller eigenvalue estimate}
        \lambda_{\min} \ll \frac{mn}{||C||_F^2}.
    \end{equation}
    Applying this estimate and \ref{Bessel estimate 1} to $J_{\ell}(4\pi\sqrt{\lambda_{\min}}\sin \theta)$, and applying the simple estimate \cite{GradshteynRyzhik2015}
    \begin{equation}\label{Bessel estimate 2}
        J_\ell(x) \ll 1, \hspace{3mm} x > 0 , \ell > \frac{1}{2}
    \end{equation}
    to $J_{\ell}(4\pi\sqrt{\lambda_{\max}}\sin \theta)$, we may cut off the sum in $R_2$ by $||C||_F \ll \frac{\sqrt{mn}}{K}$ up to an negligible error. In this range we change the summation order and deal with the inner sum
    \begin{equation}
        \sum_k g(k) J_{\ell}(4\pi\sqrt{\lambda_1}\sin \theta)J_{\ell}(4\pi\sqrt{\lambda_2}\sin \theta)
    \end{equation}
    by making use of the following integral representation of product of two Bessel functions (See equation (8) on page 47 of \cite{Transendental}):
    \begin{equation}\label{product formula}
    J_v(z)J_v(\zeta) = \frac{2}{\pi} \int_0^\frac{\pi}{2} \cos((z-\zeta)\cos \alpha)J_{2v}(2\sqrt{z \zeta}\sin\alpha)d \alpha, \hspace{3mm} \Re(v) > -\frac{1}{2}, z > 0, \zeta > 0.
    \end{equation}
    By choosing $v = \ell$, $z = 4\pi\sqrt{\lambda_1}\sin \theta$, $\zeta = 4\pi\sqrt{\lambda_2}\sin \theta$, and setting
    \begin{equation}
        \xi = 2\sqrt{z\zeta}\sin \alpha = \frac{8\pi\sqrt{mn}}{\sqrt{|\det C|}}\sin \theta \sin \alpha,
    \end{equation}
    we obtain
    \begin{equation}
        \sum_k g(k) J_{\ell}(z)J_{\ell}(\zeta) = \frac{2}{\pi} \int_0^\frac{\pi}{2} \cos((z-\zeta)\cos \alpha) \left(\sum_kg(k) J_{2k-3}(\xi)\right)d \alpha.
    \end{equation}
    Let $r = 2k-3$. We have $r \equiv 1 \bmod 4$, since $k$ is even. Setting $g_1(x) = g\left(\frac{x+3}{2}\right)$, we have
    \begin{equation}
        \sum_kg(k) J_{2k-3}(\xi) = \sum_{r \equiv 1 \bmod 4} g_1(r)J_r(\xi).
    \end{equation}
    From here the method of Newmann series can be applied, in view of the following integral representation of Bessel functions of integral order \cite{GradshteynRyzhik2015}:
    \begin{equation}
        J_r(x) = \int_{-\frac{1}{2}}^\frac{1}{2}e(rt)e^{-ix\sin 2\pi t} dt.
    \end{equation}
    We quote the following result (Lemma 5.8 in \cite{TopicsinClassicalAutomorphicForms}):
    \begin{equation}\label{Lemma 5.8}
        4\sum_{r \equiv 1 \bmod 4} g_1(r)J_r(\xi) = g_1(\xi) + h(\xi) + O(\xi c_3(g_1)),
    \end{equation}
    where $h(\xi)$ and $c_3(g_1)$ are defined by
    \begin{align}
        h(\xi) &= \int_0^\infty g_1(\sqrt{2\xi y})\sin (\xi+y-\frac{\pi}{4})(\pi y)^{-\frac{1}{2}} dy, \\
        c_3(g_1) &= \int_{-\infty}^\infty |\hat{g_1}(t)t^3| dt.
    \end{align}
    We refer readers to section 5.5 of \cite{TopicsinClassicalAutomorphicForms} for a proof of \ref{Lemma 5.8}. Recall that for $g_1$ we have
    \begin{equation}
        g_1^{(j)}(x) \ll_j K^{-j}
    \end{equation}
    for any $j \geq 0$. Thus by repeated partial integration we have
    \begin{align}
        h(\xi) &\ll_j (\xi K^{-2})^j, \label{estimate of h} \\
        c_3(g_1) &\ll K^{-3}. \label{estimate of c_3}
    \end{align}
    See also (5.73) and (5.74) in \cite{TopicsinClassicalAutomorphicForms}.
    
    The contribution of $g_1(\xi)$ to $R_2$ is
    \begin{equation}\label{R_2^{*,1}}
        R_2^{g_1(\xi)} = \sum_{||C||_F \ll \frac{\sqrt{mn}}{K}}\frac{K(nI,mI;C)}{|\det C|^\frac{3}{2}} \int_0^\frac{\pi}{2} \int_0^\frac{\pi}{2} \cos((z-\zeta)\cos \alpha) g_1(\xi) d \alpha \sin \theta d \theta.
    \end{equation}
    In view of the support of $g_1$, the sum in \ref{R_2^{*,1}} is confined in the range
    \begin{equation}
        \xi = \frac{8\pi\sqrt{mn}}{\sqrt{|\det C|}}\sin \theta \sin \alpha \gg K.
    \end{equation}
    Thus we have $\det(C) \ll \frac{mn}{K^2}$. By the estimate \ref{symplectic Kloosterman sum trivial bound}, we obtain
    \begin{equation}
        R_2^{g_1(\xi)} \ll \sum_{\substack{0 \neq |\det C| \ll \frac{mn}{K^2} \\ ||C||_F \ll \frac{\sqrt{mn}}{K}}} 1 =\sum_{0 \neq |d| \ll \frac{mn}{K^2}} \sum_{\substack{\det C = d \\ ||C|| \ll \frac{\sqrt{mn}}{K}}} 1 = \sum_{0 \neq |d| \ll \frac{mn}{K^2}} P_d\left(C \cdot \frac{mn}{K^2}\right),
    \end{equation}
    where $P_d(X)$ is the hyperbolic lattice counting function
    \begin{equation}
        P_d(X) = |\{(\alpha,\beta,\gamma,\delta) \in \bZ^4: \alpha \delta - \beta \gamma = d, \alpha^2 + \beta^2 + \gamma^2 + \delta^2 \leq X\}|,
    \end{equation}
    and $C>0$ is some constant. For $1 \leq d \leq X$ we have the following asymptotic formula (see Theorem 12.4 in \cite{SpectralMethods}):
    \begin{equation}
        P_d(X) = 6\left(\sum_{\tau|d}\tau^{-1}\right)(X + O(d^\frac{1}{3}X^\frac{2}{3})) \ll X \log |d|.
    \end{equation}
    This estimate also applies to $-X \leq d \leq -1$ by symmetry. Thus we have
    \begin{equation}
        R_2^{g_1(\xi)} \leq \frac{mn}{K^2} \sum_{0 \neq |d| \ll \frac{mn}{K^2}} \log |d| \ll \frac{m^2n^2}{K^4}\log \frac{mn}{K^2} \ll_\epsilon \left(\frac{(mn)^{2+\epsilon}}{K^{4+2\epsilon}}\right).
    \end{equation}
    
    The contributions of $h(\xi)$ and $O(\xi c_3(g_1))$ are analyzed similarly, making use the bounds \ref{estimate of h} and $\ref{estimate of c_3}$. We have
    \begin{align}
        R_2^{h(\xi)} &\ll_{j} \frac{(mn)^{\frac{j}{2}+1}}{K^{2j+2}}, \\
        R_2^{O(\xi c_3(g_1))} &\ll_\epsilon \frac{(mn)^{2+\epsilon}}{K^{6+2\epsilon}},
    \end{align}
    for any $j \geq 3$ and small $\epsilon>0$. Thus we obtain
    \begin{equation}
        R_2 \ll _{j,\epsilon} \frac{(mn)^{2+\epsilon}}{K^{4+2\epsilon}} + \frac{(mn)^{\frac{j}{2}+1}}{K^{2j+2}}.
    \end{equation}
    
    Combining the estimates of $R_1$ and $R_2$ above, and that 
    \begin{equation}
        \sum_k g(k) = \sum_k \Omega\left(\frac{k}{K}\right) \gg K,
    \end{equation}
    by our choice of $\Omega$, the proof is now complete.
\end{proof}

\section{Proof of Main Theorems}\label{section 4}
In this section we prove Theorems \ref{spin L-function}-\ref{average standard L-function}. We assume $\hat{\Phi}$ is supported in $(-\alpha,\alpha)$. We also set $\ell = k-3/2$ to save notations.
\subsection{Proof of Theorem 1.1}
By a standard explicit formula argument (see, for example, \cite{IwaniecLuoSarnak2000}), we can write the density function $D(F;\Phi;\text{spin})$ as 
\begin{equation}
    D(F;\Phi;\text{spin}) = \frac{2}{2\pi i}\int_{(2)}\Phi\left(\frac{s-\frac{1}{2}}{2\pi i}\log k^2\right)\frac{\Lambda^\prime}{\Lambda}(s,F;\text{spin})ds + 2\Phi\left(\frac{\log k^2}{2\pi i}\right)\delta_{F \in H_k^*(\Gamma_2)},
\end{equation}
where $\delta_{F \in H_k^*(\Gamma_2)} = 1$ if $F \in H_k^*(\Gamma_2)$ is a Saito-Kurokawa lift (in which case $L(s,F;\text{spin})$ has a pole at $s=3/2$) and is 0 otherwise. By \ref{spin L-function functional equation} and \ref{spinor L-function Euler product} we may further write
\begin{equation}
\begin{split}
    D(F;\Phi;\text{spin}) &= \frac{2}{\log k^2}\int_\bR \Phi(x)\left(-\log (2\pi)^2 + \frac{\Gamma^\prime}{\Gamma}\left(1+\frac{2\pi ix}{\log k^2}\right) + \frac{\Gamma^\prime}{\Gamma}\left(k-1+\frac{2\pi ix}{\log k^2}\right)\right)dx \\
    &-\frac{2}{\log k^2}\sum_{m=1}^\infty\sum_p c_m(p;F)\frac{\log p}{p^{m/2}}\hat{\Phi}\left(\frac{m\log p}{\log k^2}\right) \\
    &+2\Phi\left(\frac{\log k^2}{2\pi i}\right)\delta_{F \in H_k^*(\Gamma_2)}
\end{split}
\end{equation}
by shifting contour from $\sigma = 2$ to $\sigma = 1/2$.

For the integral involving gamma factors, we use the following estimates \cite{GradshteynRyzhik2015}:
\begin{align}
    \frac{\Gamma^\prime}{\Gamma}(a+bi) + \frac{\Gamma^\prime}{\Gamma}(a-bi) &= 2\frac{\Gamma^\prime}{\Gamma}(a) + O\left(\frac{b^2}{a^2}\right), \hspace{3mm} a > 0, b \in \bR, \\
    \frac{\Gamma^\prime}{\Gamma}(k-1) &= \log k + O\left(\frac{1}{k}\right)
\end{align}
to get
\begin{equation}
    \frac{2}{\log k^2}\int_\bR \Phi(x)\left(-\log (2\pi)^2 + \frac{\Gamma^\prime}{\Gamma}\left(1+\frac{2\pi ix}{\log k^2}\right) + \frac{\Gamma^\prime}{\Gamma}\left(k-1+\frac{2\pi ix}{\log k^2}\right)\right)dx = \hat{\Phi}(0) + o(1).
\end{equation}

For $c_1(p;F)$ and $c_2(p;F)$ we sum over $F$ against the weight $\omega_F$. Using \ref{harmonic weight}, Lemma \ref{combinatorial lemma} and the observation \ref{observation} we obtain
\begin{align}
    \sum_{F \in H_k(\Gamma_2)} \omega_F c_1(p;F) &= \Delta_k(pI,I) + \frac{\lambda_p}{\sqrt{p}}\Delta_k(I,I), \\
    \sum_{F \in H_k(\Gamma_2)} \omega_F c_2(p;F) &= -\Delta_k(pI,pI) + 2\Delta_k(p^2I,I) + O\left(\frac{1}{p}\right)\Delta_k(I,I).
\end{align}

Collecting these we have the "explicit formula":
\begin{equation}\label{sum for spinor L-functions}
\begin{split}
    \sum_{F \in H_k(\Gamma_2)}\omega_F D(F;\Phi;\text{spin}) &= \hat{\Phi}(0) + o(1) \\
    &-\frac{2}{\log k^2}\sum_p\left(\Delta_k(pI,I) + \frac{\lambda_p}{\sqrt{p}}\Delta_k(I,I)\right)\frac{\log p}{\sqrt{p}}\hat{\Phi}\left(\frac{\log p}{\log k^2}\right) \\
    &-\frac{2}{\log k^2}\sum_p\left(-\Delta_k(pI,pI) + 2\Delta_k(p^2I,I) + O\left(\frac{1}{p}\right)\Delta_k(I,I)\right)\frac{\log p}{p}\hat{\Phi}\left(\frac{2\log p}{\log k^2}\right) \\
    &-\frac{2}{\log k^2}\sum_{m=3}^\infty\sum_p \left(\sum_{F \in H_k(\Gamma_2)}\omega_F c_m(p;F)\right)\frac{\log p}{p^{m/2}}\hat{\Phi}\left(\frac{m\log p}{\log k^2}\right) \\
    &+2\Phi\left(\frac{\log k^2}{2\pi i}\right)\sum_{F_f \in H_k^*(\Gamma_2)}\omega_{F_f}
\end{split}
\end{equation}

We treat the terms $\Delta_k(pI,I), \Delta_k(pI,pI)$ and $\Delta_k(p^2I,I)$ using the Petersson formula \ref{Petersson formula}. Take the term $\Delta_k(pI,I)$ for example:
\begin{equation}
    \Delta_k(pI,I) = \frac{\sqrt{2}\pi}{8}G_{1,k}(pI,I) + \pi^2G_{2,k}(pI,I).
\end{equation}

The rank 1 term $G_{1,k}(pI,I)$ is 
\begin{equation}
    G_{1,k}(pI,I) = \sum_{\pm}\sum_{s=1}^\infty\sum_{c=1}^\infty\sum_{U,V}\frac{(-1)^{\frac{k}{2}}}{c^\frac{3}{2}(ps)^\frac{1}{2}}H^{\pm}(UU^T,pV^{-1}V^{-T};c)J_{\ell}\left(\frac{4\pi}{cs}\right),
\end{equation}
after a change of variable $s \mapsto ps$, where the summation $\sum_{U,V}$ is over
\begin{equation}
    u_{21}^2 + u_{22}^2 = ps, \hspace{3mm} v_{11}^2 + v_{21}^2 = s.
\end{equation}
By the estimates \ref{trivial estimate H}, \ref{Bessel estimate 1} and \ref{number of integral solutions}, we bound $G_{1,k}(pI,I)$ as
\begin{equation}
    G_{1,k}(pI,I) \ll \sum_{s=1}^\infty \sum_{c=1}^\infty (ps)^\epsilon s^\epsilon c^{-\frac{3}{2}}p^{-\frac{1}{2}}s^{-\frac{1}{2}}c^2\left(\frac{4\pi}{cs\ell}\right)^{\ell} \ll p^{-\frac{1}{2}+\epsilon}\left(\frac{4\pi}{\ell}\right)^{\ell},
\end{equation}
for $k$ sufficiently large. Thus its contribution to \ref{sum for spinor L-functions} is at most
\begin{equation}
\begin{split}
    \frac{1}{\log k}\sum_p p^{-\frac{1}{2}+\epsilon} \left(\frac{4\pi}{\ell}\right)^\ell \frac{\log p}{\sqrt{p}}\hat{\Phi}\left(\frac{\log p}{\log k^2}\right) &\ll \frac{1}{\log k} \left(\frac{4\pi}{\ell}\right)^\ell \sum_{p \leq k^{2\alpha}} p^{-1+\epsilon} \log p \\
    &\ll \frac{1}{\log k} \left(\frac{4\pi}{\ell}\right)^\ell k^{2\alpha \epsilon} = o(1)
\end{split}
\end{equation}
for any $\alpha > 0$ as $k \to \infty$.

The rank 2 term $G_{2,k}(pI,I)$ is 
\begin{equation}
    G_{2,k}(pI,I) = \sum_{C \in \cC} \frac{K(I,pI;C)}{|\det C|^{3/2}} \int_0^\frac{\pi}{2} J_\ell(4\pi \sqrt{\lambda_{\min}} \sin \theta) J_\ell(4\pi \sqrt{\lambda_{\max}} \sin \theta)\sin \theta d \theta,
\end{equation}
where $\lambda_{\min}$ and $\lambda_{\max}$ are the smaller and larger eigenvalues of $pC^{-1}C^{-T}$ respectively. By estimates \ref{symplectic Kloosterman sum trivial bound}, \ref{Bessel estimate 1}, \ref{Bessel estimate 2} and \ref{smaller eigenvalue estimate}, we have
\begin{equation}
    G_{2,k}(pI,I) \ll \sum_{C \in \cC} \left(\frac{4\pi\sqrt{p}\sin \theta}{\ell ||C||_F}\right)^\ell \ll p^{\frac{\ell}{2}} \left(\frac{4\pi}{\ell}\right)^\ell,
\end{equation}
for $k$ sufficiently large. Thus its contribution to \ref{sum for spinor L-functions} is at most 
\begin{equation}
    \frac{1}{\log k}\sum_p p^{\frac{\ell}{2}} \left(\frac{4\pi}{\ell}\right)^\ell \frac{\log p}{\sqrt{p}}\hat{\Phi}\left(\frac{\log p}{\log k^2}\right) \ll \frac{1}{\log k} \left(\frac{4\pi}{\ell}\right)^\ell \sum_{p \leq k^{2\alpha}} p^{\frac{\ell-1}{2}} \log p \ll k^\alpha \left(\frac{4\pi k^\alpha}{\ell}\right)^\ell,
\end{equation}
which goes to 0 as $k \to \infty$ when $\alpha < 1$. Thus we have proved the contribution of $\Delta_k(I,pI)$ to \ref{sum for spinor L-functions} is small when $\alpha < 1$. Other off-diagonal contributions are estimated similarly, and are all small when $\alpha < 1$. We skip the details here. 

The diagonal contribution of $\frac{\lambda_p}{\sqrt{p}}\Delta_k(I,I)$ to \ref{sum for spinor L-functions} from the $m=1$ term is
\begin{equation}
\begin{split}
    -\frac{2}{\log k^2}\sum_p \frac{\lambda_p \log p}{p} \hat{\Phi}\left(\frac{\log p}{\log k^2}\right) &= -\frac{2}{\log k^2}\int_1^\infty \frac{\log x}{x} \hat{\Phi}\left(\frac{\log x}{\log k^2}\right) d\pi(x) \\
    &= -\frac{2}{\log k^2}\int_1^\infty \frac{\log x}{x} \hat{\Phi}\left(\frac{\log x}{\log k^2}\right) \frac{1}{\log x} dx + o(1) \\
    &= -2\int_0^\infty \hat{\Phi}(y) dy + o(1) \\
    &= -\Phi(0) + o(1).
\end{split}
\end{equation}
Here we have used the Prime Number Theorem (PNT) for the prime counting function $\pi(x)$ and the fact the $\lambda_p = 1 + \chi_{-4}(p)$ takes values 0 and 2 for primes $p$ with density $\frac{1}{2}$ each. 

The diagonal contribution of $-\Delta_k(pI,pI)$ from the $m=2$ term is
\begin{equation}
    \frac{2}{\log k^2} \sum_p \frac{\log p}{p} \hat{\Phi}\left(\frac{2\log p}{\log k^2}\right) = \frac{\hat{\Phi}(0)}{2} + o(1).
\end{equation}

It is shown in \cite{KowalksiSahaTsimerman2012} that the $m \geq 3$ terms in \ref{sum for spinor L-functions} contribute at most $O\left(\frac{1}{\log k}\right)$. For non-Saito-Kurokawa lifts $F$ this follows from the Ramanujan bound $|c_m(p;F)| \leq 4$ and the fact
\begin{equation}
    \sum_{m=3}^\infty \sum_p \frac{\log p}{p^{m/2}} < \infty. 
\end{equation}
For a treatment of Saito-Kurokawa lifts, we refer readers to section 5 of \cite{KowalksiSahaTsimerman2012}. 

For the last term in \ref{sum for spinor L-functions} we use the fact that (see for example page 1754 of \cite{Blomer2019}) 
\begin{equation}
    \omega_{F_f} \ll \frac{1}{k^3} \frac{L(1/2,f \times \chi_{-4})}{L(1,\text{sym}^2 f)}.
\end{equation}
This, combined with the convexity bound for $L(1/2,f \times \chi_{-4})$ and the lower bound \cite{HoffsteinLockhart1994}
\begin{equation}
    L(1,\text{sym}^2 f) \gg k^{-\epsilon}
\end{equation}
give us
\begin{equation}\label{Saito-Kurokawa lift contribution}
    \sum_{F_f \in H_k^*(\Gamma_2)}\omega_{F_f} = o(1).
\end{equation}

Combining all results above, we finally have
\begin{equation}
    \sum_{F \in H_k(\Gamma_2)}\omega_F D(F;\Phi;\text{spin}) = \hat{\Phi}(0) - \Phi(0) + \frac{\Phi(0)}{2} + o(1) = \hat{\Phi}(0) - \frac{\Phi(0)}{2} + o(1).
\end{equation}
for $\alpha < 1$, as $k \to \infty$. This completes the proof of Theorem \ref{spin L-function}.

\subsection{Proof of Theorem 1.2}
The proof is similar to that of Theorem \ref{spin L-function}. We have seen in the proof there that the symmetry type is determined by diagonal contributions. So we shall concentrate on those terms and be brief about the rest. 

By Lemma \ref{combinatorial lemma}, the $m=1$ term is
\begin{equation}
    \sum_{F \in H_k(\Gamma_2)}\omega_F\tau_2(p;F) = \Delta_{k}(pI,pI) - \Delta_k(p^2I,I) + \frac{\lambda_p}{\sqrt{p}}\Delta_k(pI,I) + \left(\frac{\mu_p}{p}-1\right)\Delta_k(I,I),
\end{equation}
in which the diagonal contribution of $\Delta_k(pI,pI)$ and $-\Delta_k(I,I)$ cancel each other. 

The $m=2$ term is
\begin{equation}
\begin{split}
    \sum_{F \in H_k(\Gamma_2)}\omega_F\tau_4(p;F) &= -\Delta_k(p^3I,pI) + \Delta_k(p^2I,p^2I) + \frac{\lambda_p}{\sqrt{p}}\Delta_k(p^2I,pI) + \left(\frac{\mu_p}{p}-1\right)\Delta_k(pI,pI) \\
    &-\frac{\lambda_p}{\sqrt{p}}\Delta_k(p^3I,I) + O\left(\frac{1}{p}\right)\Delta_k(p^2I,I) + O\left(\frac{1}{\sqrt{p}}\right)\Delta_k(pI,I) + \left(O\left(\frac{1}{p}\right)+1\right)\Delta_k(I,I).
\end{split}
\end{equation}
The diagonal contribution from $\Delta_k(p^2I,p^2I), -\Delta_k(pI,pI)$ and $\Delta_k(I,I)$ is
\begin{equation}
    -\frac{2}{\log k^4} \sum_p \frac{\log p}{p} \hat{\Phi}\left(\frac{2\log p}{\log k^4}\right) = -\frac{\hat{\Phi}(0)}{2} + o(1).
\end{equation}

To illustrate why the range of support is restricted to $(-\frac{1}{4},\frac{1}{4})$, we analyze the contribution of the rank 2 term $G_{2,k}(pI,pI)$:
\begin{equation}
    G_{2,k}(pI,pI) = \sum_{C \in \cC} \frac{K(pI,pI;C)}{|\det C|^{3/2}} \int_0^\frac{\pi}{2} J_\ell(4\pi \sqrt{\lambda_{\min}} \sin \theta) J_\ell(4\pi \sqrt{\lambda_{\max}} \sin \theta)\sin \theta d \theta,
\end{equation}
where $\lambda_{\min}$ and $\lambda_{\max}$ are eigenvalues of $p^2C^{-1}C^{-T}$. We estimate as before to get
\begin{equation}
    G_{2,k}(pI,pI) \ll \sum_{C \in \cC} \left(\frac{4\pi p\sin \theta}{\ell ||C||_F}\right)^\ell \ll p^\ell \left(\frac{4\pi}{\ell}\right)^\ell.
\end{equation}
Thus it contributes at most
\begin{equation}
    \frac{1}{\log k}\sum_p p^\ell \left(\frac{4\pi}{\ell}\right)^\ell \frac{\log p}{\sqrt{p}}\hat{\Phi}\left(\frac{\log p}{\log k^4}\right) \ll \frac{1}{\log k} \left(\frac{4\pi}{\ell}\right)^\ell \sum_{p \leq k^{4\alpha}} p^{\ell-\frac{1}{2}} \log p \ll k^{2\alpha} \left(\frac{4\pi k^{4\alpha}}{\ell}\right)^\ell,
\end{equation}
which is $o(1)$ as $k \to \infty$ if $\alpha < \frac{1}{4}$. Other off-diagonal terms are estimated similarly. 

\subsection{Proof of Theorem 1.3}
The contribution of gamma factors and the diagonal contribution do not change upon averaging over $k$ with respect to $\Omega$. To illustrate how we may extend the range of support form $(-\frac{1}{4},\frac{1}{4})$ to $(-\frac{5}{18},\frac{5}{18})$, we take the term $\Delta_k(pI,pI)$ for example. 

By Lemma \ref{averaged Petersson}, the off-diagonal part of 
\begin{equation}
    \left(\sum_k \Omega\left(\frac{k}{K}\right)\right)^{-1} \sum_k \Omega\left(\frac{k}{K}\right)\Delta_k(pI,pI)
\end{equation}
is at most
\begin{equation}
    \frac{p}{K^4} + \frac{p^{4+2\epsilon}}{K^{5+2\epsilon}} + \frac{p^{j+2}}{K^{2j+3}}.
\end{equation}
for any $j \geq 3$ and small $\epsilon > 0$. It contributes at most 
\begin{equation}
    \frac{1}{\log K}\sum_{p \leq K^{4\alpha}} \left(\frac{p}{K^4} + \frac{p^{4+2\epsilon}}{K^{5+2\epsilon}} + \frac{p^{j+2}}{K^{2j+3}}\right)\frac{\log p}{\sqrt{p}} \ll K^{6\alpha - 4} + K^{18\alpha -5 +6\epsilon} + K^{(4j+10)\alpha - (2j+3)},
\end{equation}
which is $o(1)$ if $\alpha < \frac{5}{18}$, by taking $j$ sufficiently large. 

Finally, to see 
\begin{equation}
    \hat{\Phi}(0) - \frac{\Phi(0)}{2} = \int_{-\infty}^\infty \Phi(x)W(Sp)(x)dx,
\end{equation}
we use the Plancheral formula
\begin{equation}
    \int_{-\infty}^\infty \Phi(x)W(Sp)(x)dx = \int_{-\infty}^\infty \hat{\Phi}(y)\hat{W}(Sp)(y)dy
\end{equation}
and the Fourier pair \ref{Fourier pair}. The proof is now complete.

\section{Application to Non-vanishing of Central Values}\label{section 5}
\subsection{Proof of Corollary 1.1} 
By Theorem \ref{spin L-function} we have
\begin{equation}
    \sum_{F \in H_k(\Gamma_2)} \omega_F \sum_{\rho_{F,\text{spin}}} \Phi\left(\frac{\gamma_{F,\text{spin}}}{2\pi}\log k^2\right) < \hat{\Phi}(0) - \frac{\Phi(0)}{2} + \epsilon,
\end{equation}
for any $\epsilon>0$ and $k$ large enough. We further assume 
\begin{equation}
    \Phi(x) \geq 0, \hspace{3mm} \Phi(0) = 1.
\end{equation}
By these conditions we may pick up only the zeros $\rho_{F,\text{spin}} = 1/2$ to get
\begin{equation}
    \begin{split}
        \sum_{F \in H_k(\Gamma_2)} \omega_F \sum_{\rho_{F,\text{spin}}} \Phi\left(\frac{\gamma_{F,\text{spin}}}{2\pi}\log k^2\right) &\geq \sum_{F \in H_k(\Gamma_2)} \omega_F \cdot \text{ord}_{s=1/2} L(s,F;\text{spin}) \\
        &= \sum_{m=2}^\infty m \sum_{\text{ord}_{s=1/2} L(s,F;\text{spin}) = m} \omega_F \\
        &\geq 2\sum_{\text{ord}_{s=1/2} L(s,F;\text{spin}) \geq 2} \omega_F.
    \end{split}
\end{equation}
Here we have used the fact that the root number of $L(s,F;\text{spin})$ is always $+1$. Thus the vanishing order of $L(s,F;\text{spin})$ at $s=1/2$ is even. These inequalities, together with \ref{0th moment of harmonic weight}, give us
\begin{equation}
    \sum_{L(1/2,F;\text{spin}) \neq 0} \omega_F > 1 - \frac{1}{2}\left(\hat{\Phi}(0) - \frac{\Phi(0)}{2}\right) - \epsilon.
\end{equation}

It is discussed in \cite{IwaniecLuoSarnak2000} that the Fourier pair
\begin{equation}
    \Phi(x) = \left(\frac{\sin \pi vx}{\pi vx}\right)^2; \hspace{3mm} \hat{\Phi}(y) = \frac{1}{v}\left(1-\frac{|y|}{v}\right), \hspace{3mm} |y| < v \hspace{3mm} (v > 0).
\end{equation}
gives essentially the optimal bound. With this choice we have
\begin{equation}
    \sum_{L(1/2,F;\text{spin}) \neq 0} \omega_F > \frac{5}{4} - \frac{1}{2v} - \epsilon,
\end{equation}
for any $0 < v < 1$. Taking liminf in $k$ and $v \to 1$, we have
\begin{equation}
    \liminf_{k \to \infty} \sum_{L(1/2,F;\text{spin}) \neq 0} \omega_F \geq \frac{3}{4}.
\end{equation}
We can further ignore the contribution of Saito-Kurokawa lifts, in view of \ref{Saito-Kurokawa lift contribution}. This completes the proof of Corollary \ref{non-vanishing}.

\subsection{Further Discussion}
From the proof of Corollary \ref{non-vanishing} we see that in order to obtain any result on non-vanishing of central values of $L(s,F;\text{spin})$ or $L(s,F;\text{std})$, the range of support in the corresponding Density Theorem must go beyond $(-\frac{2}{5},\frac{2}{5})$. This range is by setting
\begin{equation}
    \frac{5}{4} - \frac{1}{2v} = 0.
\end{equation}

The previous range of support $(-\frac{4}{15},\frac{4}{15})$ obtained in \cite{KowalksiSahaTsimerman2012} for spinor $L$-functions is not large enough, for $\frac{4}{15} < \frac{2}{5}$. Thus our extension to $(-1,1)$ is significant for the purpose of non-vanishing. 

Unfortunately, for standard $L$-functions, our range of support is still not large enough to obtain a non-vanishing result, even after performing an average over weight ($\frac{5}{18} < \frac{2}{5}$). The author would like to address this problem by establishing a more refined version of Lemma \ref{averaged Petersson} in the future. 

\section*{Acknowledgement}
The author would like to thank Prof. Wenzhi Luo for his suggestion on the topic and his valuable comments.

\printbibliography

\vspace{5mm}

\noindent Department of Mathematics, The Ohio State University, Columbus, OH 43210, USA.

\noindent E-mail address: {\tt zhao.3326@buckeyemail.osu.edu}

\end{document}